\documentclass[11pt]{amsart}
\usepackage{amsmath}
\usepackage{amssymb}
\usepackage{amsthm}
\usepackage{epsfig}
\usepackage{tabls}
\usepackage{graphicx}
\usepackage{color}
\usepackage{epstopdf}
\DeclareGraphicsExtensions{.pdf,.eps}

\numberwithin{equation}{section}
\numberwithin{figure}{section}

\newtheorem{theorem}{Theorem}
\newtheorem{lemma}{Lemma}
\newtheorem{corollary}{Corollary}
\newtheorem*{zaslavsky}{Zaslavsky's Theorem}
\newtheorem*{shitheorem}{Shi's Theorem}
\theoremstyle{definition}
\newtheorem{definition}{Definition}
\newtheorem*{problem}{Open Problem}
\newtheorem*{assertions}{Assertions}

\newtheorem*{status}{Status}
\newtheorem{conjecture}{Conjecture}
\numberwithin{definition}{section}

\newcommand{\R}{\mathbb{R}}
\newcommand{\Q}{\mathbb{Q}}
\renewcommand{\DH}{\mathcal{DH}}

\renewcommand{\S}{\mathfrak{S}}
\newcommand{\I}{\mathcal{I}}
\newcommand{\A}{\mathcal{A}}
\newcommand{\F}{\mathbb{F}}
\newcommand{\Z}{\mathbb{Z}}
\newcommand{\Shi}{{\sf Shi}}

\newcommand{\Ish}{{\sf Ish}}
\newcommand{\shi}{{\sf shi}}
\newcommand{\ish}{{\sf ish}}
\newcommand{\stat}{{\sf stat}}

\newcommand{\rem}{{\sf rem}}
\newcommand{\quo}{{\sf quo}}
\newcommand{\area}{{\sf area}}
\newcommand{\bounce}{{\sf bounce}}
\newcommand{\Cox}{{\sf Cox}}
\newcommand{\Aff}{{\sf Aff}}
\newcommand{\w}{\tilde{w}}
\renewcommand{\u}{\tilde{u}}

\title{Hyperplane Arrangements and Diagonal Harmonics}
\thanks{The author learned in January 2010 that a version of the main result (Theorem \ref{th:main}) was known to Mark Haiman in 2007, however it was not published. We anticipate that the two projects will be merged in a forthcoming joint paper.}
\author{Drew Armstrong}
\address{Department of Mathematics, University of Miami, Coral Gables, FL 33146}
\email{armstrong@math.miami.edu}
\date{May 2010} 
\keywords{Shi arrangement, Ish arrangement, affine permutations, diagonal harmonics, Catalan numbers, nabla operator, parking functions}
\subjclass[2000]{05E10, 52C35}

\begin{document}
\maketitle

\begin{abstract}
In 2003, Haglund's {\sf bounce} statistic gave the first combinatorial interpretation of the $q,t$-Catalan numbers and the Hilbert series of diagonal harmonics. In this paper we propose a new combinatorial interpretation in terms of the affine Weyl group of type $A$. In particular, we define two statistics on affine permutations; one in terms of the Shi hyperplane arrangement, and one in terms of a new arrangement --- which we call the Ish arrangement. We prove that our statistics are equivalent to the {\sf area'} and {\sf bounce} statistics of Haglund and Loehr. In this setting, we observe that {\sf bounce} is naturally expressed as a statistic on the root lattice. We extend our statistics in two directions: to ``extended'' Shi arrangements and to the bounded chambers of these arrangements. This leads to a (conjectural) combinatorial interpretation for all integral powers of the Bergeron-Garsia nabla operator applied to the elementary symmetric functions.
\end{abstract}

\section{Introduction}
First we define the diagonal harmonics --- which we will keep in mind throughout --- then we will discuss hyperplane arrangements --- which the paper is really about.

\subsection{Diagonal Harmonics}
The symmetric group $\S(n)$ acts on the polynomial ring $S=\Q[x_1,\ldots,x_n]$ by permuting variables. Newton showed that the subring of $\S(n)$-invariant polynomials is generated by the algebraically independent {\sf power sum polynomials}: $p_k=\sum_{i=1}^n x_i^k$ for $k=1,2,\ldots,n$. It is known that the {\sf coinvariant ring} $R=S/(p_1,\ldots,p_n)$ is a graded version of the regular representation of $\S(n)$, with Hilbert series
\begin{equation*}
\sum_{i=0}^n \dim R_i\, q^i = \prod_{j=1}^n (1+q+q^2+\cdots +q^j) = [n]_q!.
\end{equation*}
The dual ring $S^*=\Q[\partial/\partial x_1,\ldots,\partial/\partial x_n]$ acts on $S$ via the pairing $(\partial/\partial x_i)x_j=\delta_{ij}$, hence the coinvariant ring is isomorphic to the quotient $S^*/(p_1^*,\ldots,p_n^*)$, where $p_k^*=\sum_{i=1}^n (\partial/\partial x_i)^k$ for $k=1,\ldots,n$. On the other hand, this quotient  is naturally isomorphic to the submodule $H\subseteq S$ annihilated by the $p_k^*$:
\begin{equation*}
H=\{ f\in S: p_k^*\, f =0 \text{ for all } k\}.
\end{equation*}
This $H$ is called the ring of {\sf harmonic polynomials} since, in particular, $p_2^*$ is the standard Laplacian operator on $S$.

Now consider the ring $DS=\Q[x_1,\ldots,x_n,y_1,\ldots,y_n]$ of polynomials in two sets of commuting variables, together with the {\sf diagonal action} of $\S(n)$, which permutes the $x$ variables and the $y$ variables {\em simultaneously}. Weyl \cite{weyl} showed that the $\S(n)$-invariant subring of $DS$ is generated by the {\sf polarized power sums}: $p_{k,\ell}=\sum_{i=1}^n x_i^k y_i^\ell$ for all $k+\ell >0$. Hence the ring of {\sf diagonal coinvariants} $DR=DS/( p_{k,\ell}: k+\ell>0)$ is naturally isomorphic to the ring of {\sf diagonal harmonic polynomials}:
\begin{equation*}
DH =\{ f\in DS : \sum_{i=1}^n (\partial/\partial x_i)^k(\partial/\partial y_i)^\ell\, f = 0 \text{ for all } k+\ell>0\}.
\end{equation*}
The diagonal action preserves the bigrading of $DS$ by $x$-degree and $y$-degree, hence $DH$ is a bigraded $\S(n)$-module. The bigraded Hilbert series
\begin{equation}
\DH(n;q,t) := \sum_{i,j=0}^n \dim (DH)_{i,j} \,\, q^i t^j
\label{eq:qtHilbert}
\end{equation}
has beautiful and remarkable properties. The study of $\DH(n;q,t)$ was initiated by Garsia and Haiman (see \cite{Hai94}) and is today an active area of research.

\subsection{Some Arrangements}
\begin{figure}
\begin{center}
\includegraphics{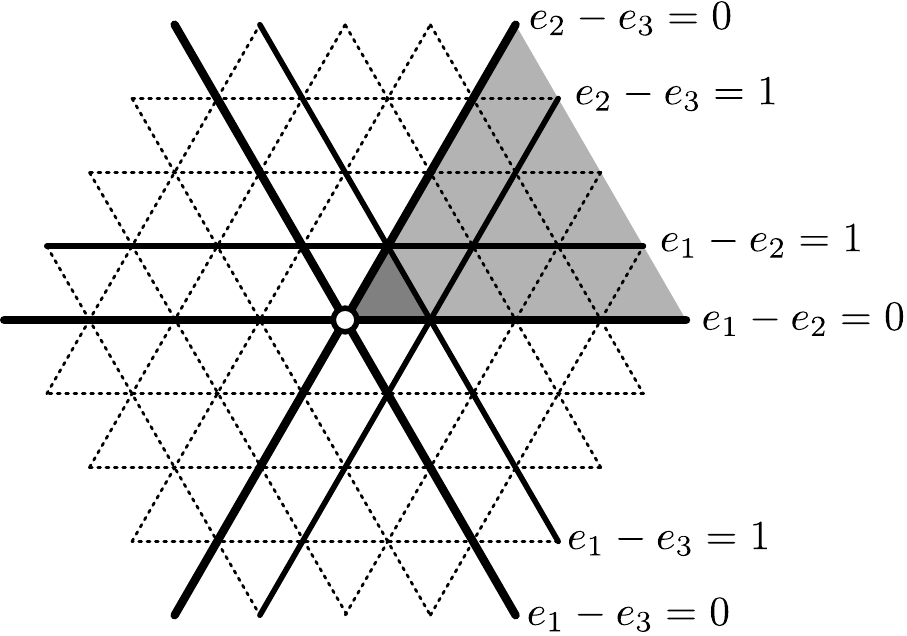}
\end{center}
\caption{Some arrangements in $\R^3_0$}
\label{fig:arrangements}
\end{figure}

Let $\{e_1,e_2,\ldots,e_n\}$ be the standard basis for $\R^n$. Given $v\in\R^n$ and $k\in\R$, we will often use the notation ``\,$v=k$\,'' as shorthand for the set $\{x: (x,v)=k\}\subseteq\R^n$, where $(\cdot,\cdot)$ is the standard inner product. Consider the following three arrangements of hyperplanes, respectively called the {\sf Coxeter arrangement}, {\sf Shi arrangement}, and {\sf affine arrangement} of type $A_{n-1}$:
\begin{align*}
\Cox(n) &:=\left\{ e_i-e_j=a : 1\leq i<j\leq n,\, a=0\right\},\\
\Shi(n)  &:=\left\{ e_i-e_j=a : 1\leq i<j\leq n,\, a\in\{0,1\}\right\},\\
\Aff(n)  &:=\left\{ e_i-e_j=a : 1\leq i<j\leq n,\, a\in\Z\right\}.
\end{align*}
Since all hyperplanes in this paper contain the line $e_1+e_2+\cdots +e_n=0$, we will typically think of these arrangements in the $(n-1)$-dimensional quotient space
\begin{equation*}
\R^n_0:= \R^n/(e_1+e_2+\cdots +e_n=0).
\end{equation*}
If $\mathcal{A}$ is an arrangement in a space $V$ then the connected components of the complement $V-\cup_{H\in\mathcal{A}} H$ are called {\sf chambers}. We will refer to chambers of the Coxeter arrangement as {\sf cones}; and refer to affine chambers as {\sf alcoves}. Let $C_\circ$ denote the {\sf dominant cone}, which satisfies the coordinate inequalities
\begin{equation*}
e_1>e_2>\cdots >e_n,
\end{equation*}
and let $A_\circ$ denote the {\sf fundamental alcove}, satisfying
\begin{equation*}
e_1>e_2>\cdots>e_n>e_1-1.
\end{equation*}
Figure \ref{fig:arrangements} displays the arrangements $\Cox(3)$, $\Shi(3)$, and $\Aff(3)$ in $\R^3_0$, with the dominant cone and fundamental alcove shaded. The Shi arrangement was introduced by Jian-Yi Shi (see \cite[Chapter 7]{shi:book}) in his description of the Kazhdan-Lusztig cells for certain affine Weyl groups.

\subsection{Symmetric Group}
The symmetric group $\S(n)$ has a faithful representation as a group of isometries of $\R^n_0$ generated by the set
\begin{equation*}
S=\{s_1,s_2,\ldots,s_{n-1}\},
\end{equation*}
where $s_i$ is the reflection in the hyperplane $e_i-e_{i+1}=0$.  The reflection $s_i$ corresponds in $\S(n)$ to the transposition of adjacent symbols $(i,i+1)$.

The symmetric group acts simply-transitively on the cones of the Coxeter arrangement $\Cox(n)$. By convention, let the dominant cone $C_\circ$ correspond to the identity permutation; then for any permutation $w\in\S(n)$ the cone $w C_\circ$ satisfies
\begin{equation*}
e_{w(1)}>e_{w(2)}>\cdots >e_{w(n)}.
\end{equation*}

\subsection{Affine Symmetric Group}
Now let $s_n$ denote the reflection in the {\em affine} hyperplane $e_1-e_n=1$. The linear reflections $\{s_1,s_2,\ldots,s_{n-1}\}$ together with the affine reflection $a_n$ generate the {\sf affine Weyl group} of type $\tilde{A}_n$. This group acts simply-transitively on the set of alcoves, where the fundamental alcove $A_\circ$ corresponds to the identity element of the group. Note that $A_\circ$ is a (non-regular) simplex in $\R^n_0$ whose facets are supported by the reflecting hyperplanes of the generators $\{s_1,s_2,\ldots,s_n\}$.

Lusztig \cite{lusztig} introduced an affine version of the symmetric group, whose  combinatorial properties were developed further by Bj\"orner and Brenti \cite{bjornerbrenti}: Let $\tilde{\S}(n)$ denote the group of infinite permutations $\w:\Z\to\Z$ satisfying:
\begin{itemize}
\item $\w(k+n)=\w(k)+n$ for all $k\in\Z$,
\item $\w(1)+\w(2)+\dots+\w(n)=\binom{n+1}{2}$.
\end{itemize}
The first property says that $\w$ is periodic and the second fixes a frame of reference. The elements of $\tilde{\S}(n)$ are called {\sf affine permutations}, and $\tilde{\S}(n)$ is the {\sf affine symmetric group}. Following Bj\"orner and Brenti, we will usually specify an affine permutation $\w\in\tilde{\S}(n)$ using the {\sf window notation}:
\begin{equation*}
``\w=[\w(1),\w(2),\cdots,\w(n)].''
\end{equation*}
For integers $i<j$ we will write $((i,j)):\Z\to\Z$ to denote the ``affine tranposition'' that swaps the elements in positions $i+kn$ and $j+kn$ for all $k\in\Z$. We could also write $((i,j))=\prod_k (i+kn,j+kn)$. Lusztig proved that the correspondence $s_i\leftrightarrow ((i,i+1))$ defines an isomorphism between the affine symmetric group and the affine Weyl group of type $\tilde{A}_n$. Here the affine tranposition $((i,j))$ corresponds to the reflection in the affine hyperplane
\begin{equation}
e_{\rem(j,n)}-e_{\rem(i,n)} = \quo(j,n),
\label{eq:affhyp}
\end{equation}
where $\quo(x,n)$ and $\rem(x,n)$ are the unique quotient and remainder of $x$ by $n$, {\bf with remainder taken in the set $\{1,2,\ldots,n\}$}. In particular, note that the generator $s_i=((i,i+1))$ corresponds to $e_i-e_{i+1}=0$ for $1\leq i\leq n-1$, and $s_n=((n,n+1))$ corresponds to $e_1-e_n=1$.

\subsection{The Ish Arrangement} 
\begin{figure}
\begin{center}
\includegraphics[scale=1]{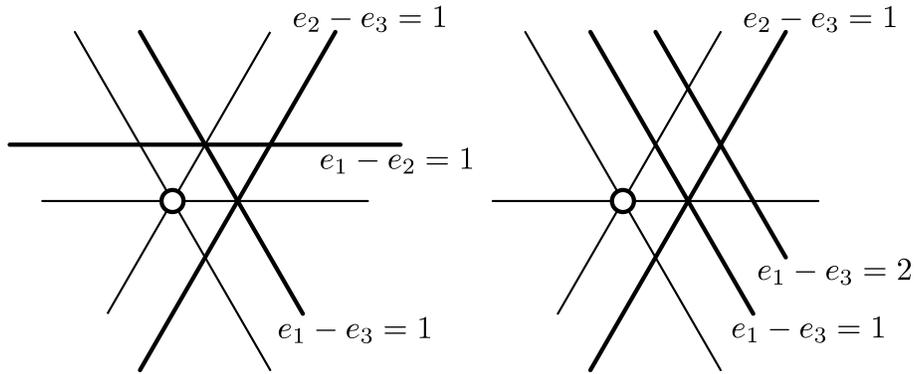}
\end{center}
\caption{The arrangements $\Shi(3)$ and $\Ish(3)$}
\label{fig:shi3ish3}
\end{figure}

To end the introduction we will introduce a new hyperplane arrangement, which we call the {\sf Ish arrangement}. Like the Shi arrangement the Ish arrangement begins with the $\binom{n}{2}$ linear hyperplanes of the Coxeter arrangement and then adds another $\binom{n}{2}$ affine hyperplanes:
\begin{equation*}
\Ish(n):=\Cox(n)\cup \{ e_i-e_n=a : 1\leq i\leq n-1,\,\, a\in\{1,\ldots,n-i\}\}.
\end{equation*}
Figure \ref{fig:shi3ish3} displays the arrangements $\Shi(3)$ and $\Ish(3)$. Note that each has $16$ chambers and $4$ bounded chambers. There is an important reason for this: the arrangements $\Shi(n)$ and $\Ish(n)$ share the same {\em characteristic polynomial}, as we now show.

To avoid extra notation, we will use a non-standard definition of the characteristic polynomial. This definition is due to Crapo and Rota, and was applied extensively by Athanasiadis  --- see Stanley \cite[Lecture 5]{stanley:hyp} for details. Let $\A$ be an arrangement of {\em finitely many} hyperplanes in $\R^n$. Suppose further that each of these hyperplanes has an equation with integer coefficients. Then, given a finite field $\F_q$ with $q$ elements, we may consider the reduced arrangement $\A_q$ in $\F_q^n$. It turns out that (for all but finitely many $q$), the number of points of $\F_q^n$ {\bf not} on any hyperplane of $\A_q$ is given by a polynomial in $q$, called the {\sf characteristic polynomial} of $\A$:
\begin{equation*}
\chi(\A,q)=\# \left( \F_q^n -\cup_{H\in\A_q} H\right) = q^n -\# \cup_{H\in\A_q} H.
\end{equation*}

The characteristic polynomial of the Shi arrangement is well known (cf. \cite[Theorem 5.16]{stanley:hyp}). Our new result is the following.
\begin{theorem}
The Shi arrangement and the Ish arrangement share the same characteristic polynomial, viz.
\begin{equation*}
\chi(\Ish(n),q)= q\,(q-n)^{n-1}.
\end{equation*}
\end{theorem}

\begin{proof}
Let $p$ be a large prime and consider a regular $p$-gon whose vertices represent the elements of the finite field $\F_p=\{1,2,\ldots,p\}$, in clockwise order. We will think of a vector $v=(v_1,\ldots,v_n)\in\F^n_p$ as a labeling of the vertices, as follows:  if $v_i=k\in\F_p$, then place the label $v_i$ on the vertex $k$.

To say that $v\in\F^n_p$ is in the complement of the reduced Ish arrangement $\Ish(n)_p$, means that $v_i-v_j\neq 0$ for all $1\leq i<j\leq n$ (that is, labels $v_i$ and $v_j$ do not occupy the same vertex) and $v_i\neq v_n+a$ for $1\leq a\leq n-i$ (that is, the label $v_i$ does not occur within the $n-i$ vertices clockwise of $v_n$). To count the vectors in the complement, first note that there are $p$ ways to place the label $v_n$. After this, we may place $v_1$ in $(p-n)$ ways, since it must avoid the position of $v_n$ and the $n-1$ positions just clockwise of this. Next, we may place $v_2$ in $(p-n)$ ways since it must avoid the position of $v_n$, the $n-2$ positions just clockwise of this, and also the position of $v_1$. Continuing in this way, we find that there are $p\, (p-n)^{n-1}$ vectors in the complement.
\end{proof}

The following is a standard result on real hyperplane arrangements.

\begin{zaslavsky}[see, e.g., Theorem 2.5 of \cite{stanley:hyp}]
Let $\A$ be an arrangement in $\R^d$ in which the intersection of all hyperplanes has dimension $k$. Then:
\begin{itemize}
\item The number of chambers of $\A$ is $(-1)^d\chi(A,-1)$.
\item The number of bounded chambers of $\A$ is $(-1)^{d-k}\chi(A,1)$.
\end{itemize}
\end{zaslavsky}
\begin{corollary}
The arrangements $\Shi(n)$ and $\Ish(n)$ have the same number of chambers --- i.e. $(n+1)^{n-1}$ --- and the same number of bounded chambers --- i.e. $(n-1)^{n-1}$.
\end{corollary}
\begin{problem}
Find a bijective proof of the corollary.
\end{problem}

The observation that the Shi and Ish arrangements are (in some undefined sense) ``dual'' to each other is at the heart of this paper.

\section{Two Statistics on Shi Chambers}
Now we define two statistics --- called $\shi$ and $\ish$ --- on the chambers of a Shi arrangement (more generally, on the elements of the group $\tilde{\S}(n)$). The first statistic is well known and the second is new. Each statistic counts a certain kind of inversions of an affine permutation, and so we begin by discussing inversions.

\subsection{Affine Inversions}
Let $w$ be an element of the (finite) symmetric group $\S(n)$. If $w(i)>w(j)$ for indices  $1\leq i<j\leq n$ we say that the tranposition $(i,j)$ is an {\sf inversion} of $w$ --- equivalently, this means that the hyperplane $e_i-e_j=0$ separates the cone $wC_\circ$ from the dominant cone $C_\circ$. The number of inversions of $w$ is called its {\sf length}.

In the affine symmetric group $\tilde{\S}(n)$, there is again a correspondence between hyperplanes and transpositions. Recall that the affine transpositions $((i,j))$ and $((i',j'))$ coincide if $i'=i+kn$ and $j'=j+kn$ for some $k\in\Z$, in which case they represent the same hyperplane \eqref{eq:affhyp}. Hence, each affine transposition has a standard representative in the set
\begin{equation*}
\tilde{T}:=\left\{ ((i,j)) :\quad 1\leq i\leq n,\quad i<j \right\} \subseteq\tilde{\S}(n).
\end{equation*}
Given an affine permutation $\w\in\tilde{\S}(n)$ and an affine transposition $((i,j))\in\tilde{T}$ such that $\w(i)>\w(j)$, we say that $((i,j))$ is an {\sf affine inversion} of $\w$ --- equivalently, the hyperplane \eqref{eq:affhyp} separates the alcove $\w A_\circ$ from the fundamental alcove $A_\circ$. Again, the {\sf (affine) length} of $\w$ is its number of affine inversions.

\subsection{The {\sf shi} statistic}
\begin{figure}
\begin{center}
\includegraphics[scale=1.3]{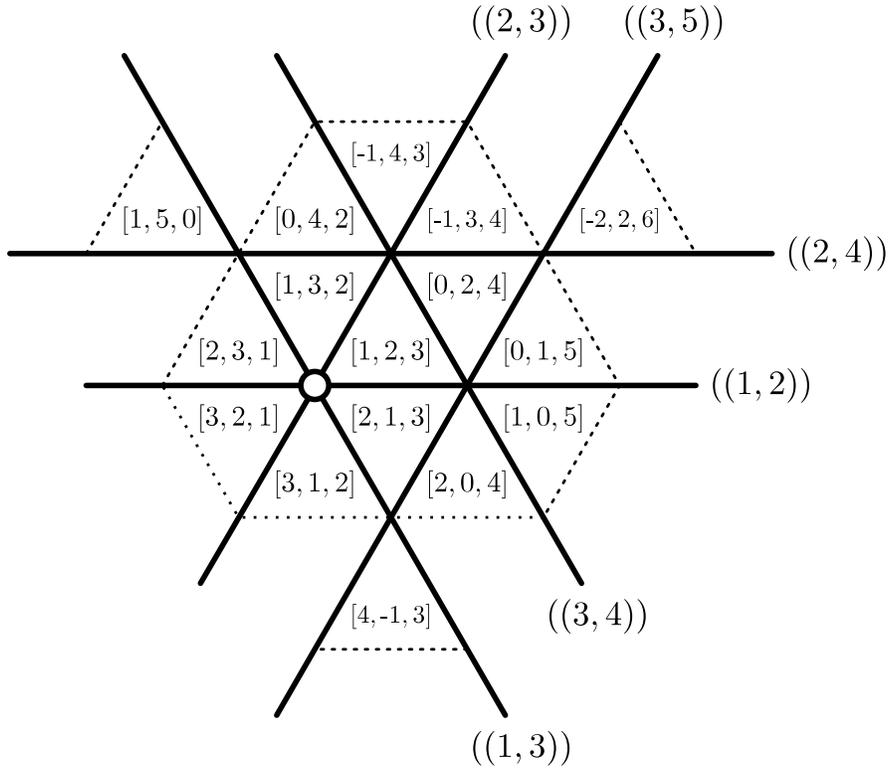}
\end{center}
\caption{Chambers of $\Shi(3)$ labeled by affine permutations}
\label{fig:affperm}
\end{figure}

Each chamber of the Shi arrangement contains a set of alcoves and we will see (Theorem \ref{th:sommers}) that among these is a unique alcove of minimum length --- which we call the {\sf representing alcove} of the chamber, or just a {\sf Shi alcove}. This defines an injection from Shi chambers into the affine symmetric group. Figure \ref{fig:affperm} displays the representing alcoves for $\Shi(3)$, labeled by affine permutations. We have labeled the Shi hyperplanes with their corresponding affine transpositions,
\begin{equation*}
\Shi(n)=\left\{ ((i,j)) : \quad1\leq i\leq n,\quad i<j<n+i\right\}.
\end{equation*}

\begin{definition}
Given a Shi chamber with representing alcove $A$, let $\shi(A)$ denote the number of Shi hyperplanes separating $A$ from the fundamental alcove $A_\circ$. Equivalently, if $A=\w A_\circ$ for affine permutation $\w\in\S(n)$, then $\shi(\w)$ is the number of affine inversions $((i,j))$ of $\w$ satisfying $i<j<n+i$.
\end{definition}

For example, consider the permutation $\w=[1,5,0]$ in the figure. The inversions of $\w$ are $((1,3)), ((2,3)), ((2,4)), ((2,6))$, and hence $\w$ has length $4$. However, only three of these --- viz. $((1,3)), ((2,3)), ((2,4))$ --- come from Shi hyperplanes, hence $\shi(\w)=3$.

\subsection{The {\sf ish} statistic}
To give a natural definition for our second statistic, we must discuss the quotient group $\tilde{\S}(n)/\S(n)$. By abuse of notation, let $\S(n)$ denote the subgroup of $\tilde{\S}(n)$ generated by the subset
\begin{equation*}
I=\{s_1,\ldots,s_{n-1}\}\subseteq\{s_1,\ldots,s_{n-1},s_n\}=S.
\end{equation*}
In the language of Coxeter groups we say that $\S(n)$ is a {\sf parabolic subgroup} of $\tilde{\S}(n)$. When $W=\tilde{\S}(n)$ the standard notation for this is to write $\S(n)=W_I$. Then each affine permutation $\w$ has a canonical decomposition
\begin{equation*}
\w=w_I\w^I,
\end{equation*}
where $w_I\in W_I$ is a finite permutation and $\w^I\in W$ is the unique coset representative of minimum (affine) length. Combinatorially, $[\w^I(1),\ldots,\w^I(n)]$ is the increasing rearrangement of $[\w(1),\ldots,\w(n)]$ and $w_I$ is the finite permutation needed to achieve the rearrangement. Geometrically, alcoves of the form $A=\w^I A_\circ$ are contained in the dominant cone $C_\circ$; hence $\w A_\circ =w_I A$ is contained in the cone $w_I C_\circ$.

We define the $\ish$  statistic in terms of minimal coset representatives.

\begin{definition}
Consider a Shi chamber with representing alcove $A$ and suppose that $A=\w A_\circ$. Its minimal coset representative $\w^I A_\circ$ is an alcove in the dominant cone $C_\circ$. Let $\ish(A)$ denote the number of hyperplanes of the form $e_i-e_n=a$  (with $1\leq i\leq n-1$ and $a\in\Z$) separating $\w^I A_\circ$ from the fundamental alcove $A_\circ$. Equivaently, let $\ish(\w)$ denote the number of affine inversions of $\w^I$ of the form $((n,j))$.
\end{definition}

Two notes: In order to facilitate later generalization, we have defined $\ish$ in terms of {\em all} hyperplanes of the form $e_i-e_n=a$. In our current context, however, only the Ish hyperplanes (i.e. $a\in\{1,\ldots,n-i\}$) will contribute. We also emphasize the fact that {\bf $\ish$ is a statistic on the (representing alcoves of) Shi chambers, not on the Ish chambers}. It seems that the chambers of the Ish arrangement are not so natural.

\begin{figure}
\begin{center}
\includegraphics[scale=1]{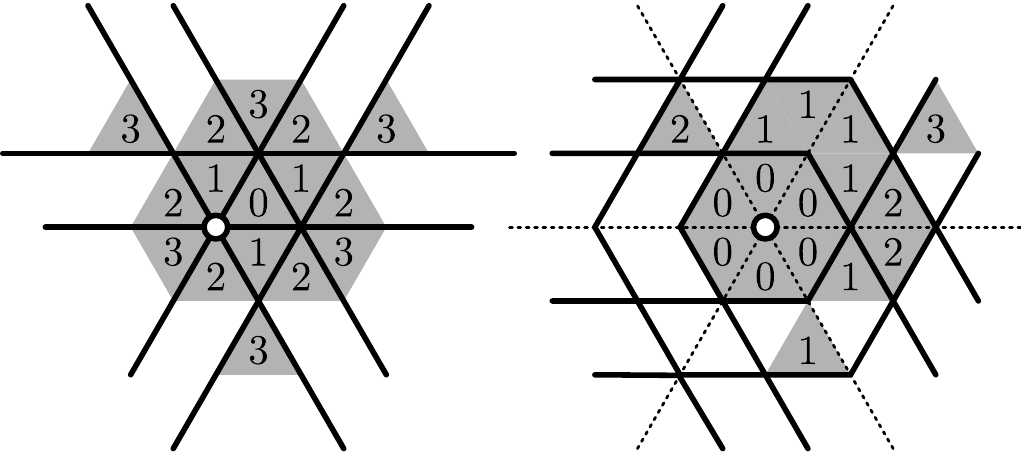}
\end{center}
\caption{The $\shi$ and $\ish$ statistics on the chambers of $\Shi(3)$}
\label{fig:shiandish}
\end{figure}

For example, consider the affine permutation $\w=[\text{-}1,4,3]$, as shown in Figure \ref{fig:affperm}. It is contained in the cone $[1,3,2]\,C_\circ$ and its increasing rearrangement is $[\text{-}1,3,4]$. Hence, it has parabolic decomposition
\begin{equation*}
[\text{-}1,4,3]=\w=w_I\w^I=[1,3,2]\,[\text{-}1,3,4].
\end{equation*}
The inversions of $\w^I=[\text{-}1,3,4]$ are $((2,4))$ and $((3,4))$, of which only the second is an Ish hyperplane; hence $\ish(\w)=1$. In Figure \ref{fig:shiandish} we have displayed the $\shi$ and $\ish$ statistics for all chambers of $\Shi(3)$. (Note: to compute $\ish$ by hand, one may extend the Ish hyperplanes from the dominant cone to the other cones by reflection.) Their joint-distribution is recorded in the following table:
\begin{center}
\begin{tabular}{rc}
& \,\,\,\,\,\,\,\,\,\,$\ish$ \\
\begin{tabular}{c} \vspace{.1in}\\ $\shi$\!\!\!\!\!\!\!\! \end{tabular} & 
\begin{tabular}{r|cccc}
 & 0 & 1 & 2 & 3\\
\hline
0 & 1& & & \\
1 & 2& 1& & \\
2 & 2& 3& 1& \\
3 & 1& 2& 2& 1
\end{tabular}\\
\end{tabular}
\end{center}

\subsection{Theorems and a Conjecture} We will make four assertions and then describe our state of knowledge about them (i.e. whether each is a Theorem or a Conjecture). We will use the following notation.

Recall from \eqref{eq:qtHilbert} that $\DH(n;q,t)$ denotes the bigraded Hilbert series of the ring of diagonal harmonic polynomials. Define
\begin{equation*}
\Shi(n;q,t):=\sum_A q^{\ish(A)} t^{\binom{n}{2}-\shi(A)},
\end{equation*}
where the sum is taken over representing alcoves $A$ for the chambers of the arrangement $\Shi(n)$. We say that an alcove is {\sf positive} if it is contained in the dominant cone $C_\circ$ (i.e. if $A$ is on the ``positive'' side of each generating hyperplane). Let $\Shi_+(n;q,t)$ denote the corresponding sum over positive Shi alcoves. Finally, consider the standard $q$-integer, $q$-factorial, and $q$-binomial coefficient:
\begin{align*}
[a]_q&=1+q+\cdots +q^{a-1},\\
[a]_q!&=[a]_q[a-1]_q\cdots [2]_q[1]_q,\\
{a \brack b}_q&=\frac{[a]_q!}{[a-b]_q![b]_q!}.
\end{align*}

\begin{assertions}\hspace{.1in}
\begin{enumerate}
\item $\Shi(n;q,t)=\DH(n;q,t)$, and hence is symmetric in $q$ and $t$.
\item $q^{\binom{n}{2}}\Shi(n;q,1/q)=[n+1]_q^{n-1}$.
\item $\Shi_+(n;q,t)$ is equal to Garsia and Haiman's $q,t$-Catalan number, and hence is symmetric in $q$ and $t$.
\item $q^{\binom{n}{2}}\Shi_+(n;q,1/q)=\frac{1}{[n]_q}{2n \brack n-1}_q$, the $q$-Catalan number.
\end{enumerate}
\end{assertions}

In particular, note that $q^{\binom{n}{2}}\Shi_+(n;q,1/q)$ is equal to the sum of $q^{\shi(A)+\ish(A)}$ over the positive Shi alcoves $A$. For $n=3$ we may compute this sum using the data in Figure \ref{fig:shiandish} to obtain
\begin{equation*}
1+q^2+q^3+q^4+q^6= \frac{[6]_q[5]_q}{[3]_q[2]_q}=\frac{1}{[3]_q}{6 \brack 2}_q,
\end{equation*}
which is a $q$-Catalan number. One may check that the other three assertions are also true in the case $n=3$.

In the following section we will establish a bijection (Theorem \ref{th:main}) from Shi chambers to labeled lattice paths, which sends our statistics $(\ish,\shi-\binom{n}{2})$ to the statistics $(\bounce,\area')$ of Haglund and Loehr \cite{haglundloehr}. This allows us to clarify the Assertions.

\begin{status} The following results all depend on our main theorem (Theorem \ref{th:main}).
\begin{enumerate}
\item {\bf Conjecture.} This is equivalent to a conjecture of Haglund and Loehr \cite{haglundloehr} (known in a different form to Haiman). No combinatorial explanation of the $q$-$t$ symmetry is known.
\item {\bf Theorem.} This is equivalent to a theorem of Loehr \cite{loehr05}.
\item {\bf Theorem.} This follows from theorems of Garsia and Haglund \cite{garsiahaglund1,garsiahaglund2}. No combinatorial explanation of the $q$-$t$ symmetry is known.
\item {\bf Theorem.} This is equivalent to a theorem of Haglund \cite{haglund03}, which was later proved bijectively by Loehr \cite{loehr07}.
\end{enumerate}
\end{status}

\section{Shi Chambers and Lattice Paths}
In this section we will prove the above stated results regarding the $\shi$ and $\ish$ statistics. To do this we will interpret Shi chambers as certain labeled lattice paths.

\subsection{The Root Poset and Dyck Paths} Cartan and Killing invented root systems prior to 1890 and used these to classify the semisimple Lie algebras. In this paper we are primarily concerned with the ``type $A$'' root system, which is related to the symmetric group. Recall that $\S(n)$ has a faithful action on $\R_0^n$ generated by the reflections $S=\{s_1,s_2,\ldots,s_{n-1}\}$, where $s_i$ is the reflection in the hyperplane $e_i-e_{i+1}=0$. The positive normal vectors to the generating hyperplanes form a special basis, called the basis of {\sf simple roots} $\Delta=\{e_1-e_2,\ldots,e_{n-1}-e_n\}$. The positive normal vectors to {\em all} reflecting hyperplanes form the set of {\sf positive roots} $\Phi^+=\{e_i-e_j : 1\leq i<j\leq n\}$.

The {\sf root poset} is a partial order on $\Phi^+$ defined as follows. Given two positive roots $\alpha,\beta\in\Phi^+$ we say that $\alpha\leq\beta$ whenever $\beta-\alpha$ can be written in the basis $\Delta$ using non-negative coefficients --- equivalently, we have $\alpha\leq\beta$ when $\beta-\alpha$ is in the positive cone generated by $\Delta$. In type $A$ this means that  $e_j-e_k\leq e_i-e_\ell$ if and only if $i\leq j<k\leq\ell$.

\begin{figure}
\begin{center}
\includegraphics[scale=1.3]{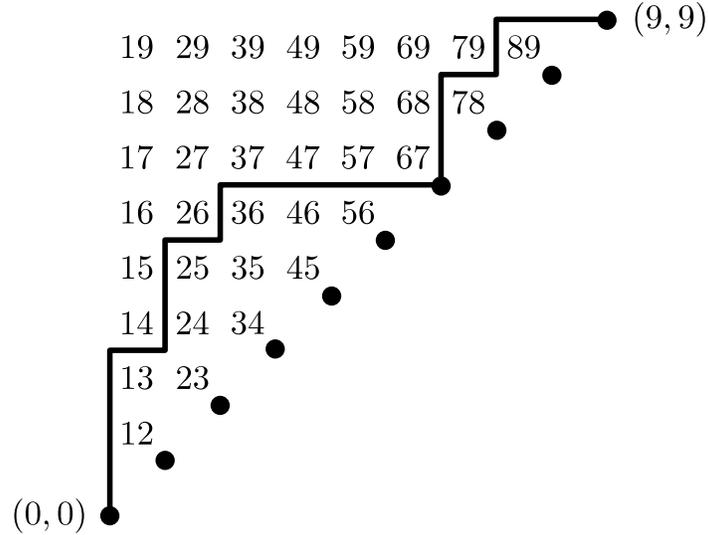}
\end{center}
\caption{An ideal and its corresponding Dyck path}
\label{fig:dyck_path}
\end{figure}

In this paper we will visualize the root poset in a particular way. Consider an array of integer points $(i,j)$, $1\leq i<j\leq n$, and place the label ``$ij$'' in the unit square with top right corner $(i,j)$. (See Figure \ref{fig:dyck_path}.) This square will represent the root $e_i-e_j$. Thus for $\alpha,\beta\in\Phi^+$ we have $\alpha\leq\beta$ when the square labeled $\beta$ occurs weakly to the left and weakly above the square labeled $\alpha$.

A set of roots $\I\subseteq\Phi^+$ is called an {\sf ideal} if $\alpha\in\I$ and $\alpha\leq\beta$ together imply $\beta\in\I$. We may picture this as a collection of unit squares aligned up and to the left. The lower boundary of these squares defines a lattice path from $(0,0)$ to $(n,n)$ which
\begin{itemize}
\item uses only steps of the form $(0,1)$ and $(1,0)$, and
\item stays weakly above the diagonal.
\end{itemize}
This defines a bijection between ideals in $\Phi^+$ and so-called {\sf Dyck paths}. For example, Figure \ref{fig:dyck_path} displays an ideal in the root poset of $\S(9)$ and its corresponding Dyck path.

\subsection{Shi Alcoves} 

\subsubsection{The Address of an Alcove}

For each root $\alpha\in\Phi^+$ and each real number $k\in\R$ let $H_{\alpha,k}$ denote the hyperplane $\{ x\in\R_0^n: (x,\alpha)=k\}$. When $\alpha=e_i-e_j$ this is the hyperplane $e_i-e_j=k$. Now let $A$ be an alcove of the affine arrangement. For each root $\alpha\in\Phi^+$ there exists a unique integer $k_A(\alpha)$ such that $A$ lies between the hyperplanes $H_{\alpha,k_A(\alpha)}$ and $H_{\alpha,k_A(\alpha)+1}$. The function $k_A:\Phi^+\to\Z$ uniquely specifies the position of $A$, so we call it the {\sf address} of $A$. An important result of J.-Y. Shi characterizes which functions can be addresses (see Sommers \cite[Proposition 4.1]{sommers}, which is a restatement of J.-Y. Shi \cite[Theorem 5.2]{shi:alcoves}).

\begin{shitheorem}
A function $k:\Phi^+\to\Z$ is the address of an alcove if and only if, for all triples $\alpha,\beta,\alpha+\beta$ of positive roots, we have
\begin{equation*}
k(\alpha)+k(\beta)\leq k(\alpha+\beta) \leq k(\alpha)+k(\beta)+1.
\end{equation*}
\end{shitheorem}	

We say that the alcove $A$ is {\sf positive} if it lies in the dominant cone $C_\circ$. Equivalently, $A$ is positive if and only if its address $k_A$ takes non-negative values.  We observe that the address of a positive alcove is an {\em increasing} function on the root poset. Indeed, if $\alpha\leq\beta$ then $\beta-\alpha$ is a non-negative integer combination of simple roots. Morever, there exists a way to get from $\alpha$ to $\beta$ by successively adding these simple roots, always staying within $\Phi^+$. Since we assumed that $k_A(\gamma)\geq 0$ for all simple $\gamma\in\Delta\subseteq\Phi^+$, the result follows from Shi's Theorem.

\subsubsection{Positive Shi Alcoves}
The Shi arrangement consists of the hyperplanes $H_{\alpha,k}$ for all $\alpha\in\Phi^+$ and $k\in\{0,1\}$. Given an alcove $A$, we would like to understand in which chamber of the Shi arrangement it occurs. This problem is easiest to solve for {\sf positive} alcoves; in this case we need only specify for which roots $k_A$ is zero and for which roots it is positive. To this end, we define
\begin{equation*}
\I_A:= k_A^{-1}(\left\{1,2,\ldots\}\right)\subseteq\Phi^+.
\end{equation*}
Since the address of a positive alcove $A$ is increasing, we observe in this case that $\I_A\subseteq\Phi^+$ is an ideal in the root poset. It turns out that this defines a bijection between positive Shi chambers and ideals. For this result we refer to Sommers \cite[Lemmas 5.1 and 5.2]{sommers}.

\begin{theorem}[Representing Alcoves]
\label{th:sommers}
Given an ideal $\I\subseteq\Phi^+$ of positive roots, there exists a unique alcove of minimum length such that $\I=\I_A$. The address of this alcove is given by $k_\I:\Phi^+\to\Z$ where $k_\I(\alpha)$ is the maximum number $r$ such that $\alpha$ can be expressed as a sum of $r$ roots in the ideal $\I$.
\end{theorem}

\begin{figure}
\begin{center}
\includegraphics[scale=1.3]{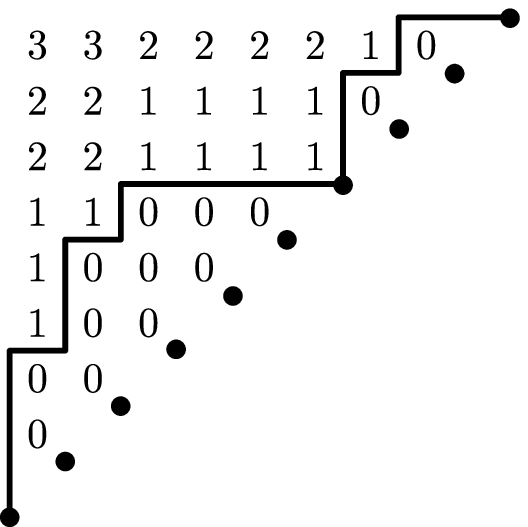}
\end{center}
\caption{The address of a positive Shi alcove}
\label{fig:alcove_address}
\end{figure}

We call the unique minimum alcove in a positive Shi chamber its {\sf representing alcove}, or just a {\sf positive Shi alcove}. Figure \ref{fig:alcove_address} displays the address of the representing alcove corresponding to the ideal in Figure \ref{fig:dyck_path}.

\subsubsection{Non-Positive Shi Alcoves}

It is true that each non-positive Shi chamber also contains a unique alcove of minimum length, which we call a {\sf non-positive Shi alcove}. Unfortunately, we do not know an expression for the address of such an alcove in the spirit of Theorem \ref{th:sommers}. Instead will use a slightly weaker result due to Pak and Stanley (see \cite[Theorem 5.1]{stanley:rota}).

Recall that a given positive Shi chamber $C$ corresponds to an ideal $\I\subseteq\Phi^+$ of positive roots: given a positive root $\alpha=e_i-e_j$, the chamber $C$ lies on the positive side of $H_{\alpha,1}$ when $\alpha\in\I$ and $C$ lies between $H_{\alpha,0}$ and $H_{\alpha,1}$ when $\alpha\not\in\I$. In addition, the {\em minimal} roots $\alpha\in\I$ (such that $\I-\alpha$ is also an ideal) correspond exactly to the hyperplanes $H_{\alpha,1}$ that support a facet of the chamber and also separate it from the fundamental alcove $A_\circ$. We call these the {\sf floors} of the chamber. In the language of Dyck paths, these  are the squares contained in the ``valleys'' of the path. For instance, the valleys in Figure \ref{fig:dyck_path} contain roots $e_1-e_4$, $e_2-e_6$, $e_6-e_7$ and $e_7-e_9$.

Now consider a non-positive Coxeter cone $wC_\circ$, with $w\in\S(n)$. The Shi hyperplanes that intersect the dominant cone $C_\circ$ have the form $e_i-e_j=1$ for $1\leq i<j\leq n$. The images of these hyperplanes in $wC_\circ$ have the form $e_{w(i)}-e_{w(j)}=1$ for $1\leq i<j\leq n$, and such a plane is actually a member of the Shi arrangement exactly when $w(i)<w(j)$. That is, the Shi planes that intersect $wC_\circ$ correspond to the {\sf non-inversions} $(i,j)$ of $w\in\S(n)$. If we then map a positive Shi alcove into $wC_\circ$, it will remain a Shi alcove if and only if its floors continue to exist. In summary, we have the following.

\begin{theorem}[Pak and Stanley \cite{stanley:rota}]
The chambers of the Shi arrangement are in bijection with pairs $(w,\I)$ where $w\in\S(n)$ is a permutation and $\I\subseteq\Phi^+$ is an ideal of positive roots (a Dyck path) such that the minimal elements of $\I$ (labels in the valleys of the path) are non-inversions of $w$.
\end{theorem}

\begin{figure}
\begin{center}
\includegraphics[scale=1.3]{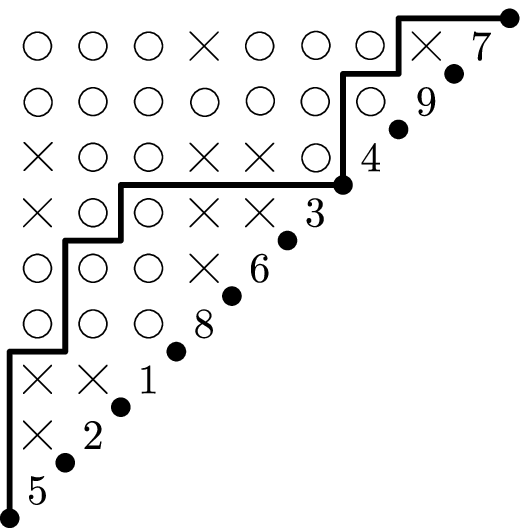}
\end{center}
\caption{A non-positive Shi chamber}
\label{fig:nonpositive_alcove}
\end{figure}

Figure \ref{fig:nonpositive_alcove} shows an example corresponding to the permutation $$w=521863497\in\S(9)$$ and the same path $\I$ as in Figures \ref{fig:dyck_path} and \ref{fig:alcove_address}. Here the symbols $\displaystyle{\times}$ and $\bigcirc$ represent, respectively, inversions and non-inversions of $w$. Notice that the valleys of $\I$ contain $\bigcirc$'s. We will call such a diagram a {\sf labeled Dyck path}.

Let us interpret the statistics $\shi$ and $\ish$ in terms of labeled Dyck paths.

\subsubsection{$\binom{n}{2}-\shi=\area'$}

In \cite{haglundloehr} Haglund and Loehr defined two statistics on labeled Dyck paths --- called $\area'$ and $\bounce$ --- and they conjectured that the generating function $\sum q^{\area'} t^{\bounce}$ equals the bigraded Hilbert series $\DH(n;q,t)$ of diagonal harmonic polynomials.

We first deal with $\area'$, which Haglund and Loehr defined as {\bf the number of non-inversions of $w$ below the labeled Dyck path $(w,\I)$}. When $w$ is the identity permutation, this is just the number of unit squares fully between the path and the diagonal, i.e. the ``area'' of the path.

\begin{theorem}
Given a Shi alcove $A$ (positive or non-positive) and its corresponding labeled Dyck path $(w,\I)$ we have
\begin{equation*}
\binom{n}{2}-\shi(A)=\area'(w,\I).
\end{equation*}
\end{theorem}

\begin{proof}
Recall that $\shi(A)$ is the number of Shi hyperplanes separating $A$ from the fundamental alcove $A_\circ$. These come in two classes. First, it is well known that the hyperplanes separating $A_\circ$ from $wA_\circ$ are exactly $e_i-e_j=0$ such that $1\leq i<j\leq n$ and $w(i)>w(j)$. These are the $\times$'s in the diagram. Second, the hyperplanes of the form $e_i-e_j=1$ separating $A$ from $wA_\circ$ correspond to unit squares above the path. Such a hyperplane is a Shi hyperplane whenever $w(i)<w(j)$, so these correspond to $\bigcirc$'s above the path. Since the total number of symbols is $\binom{n}{2}$ we conclude that $\binom{n}{2}-\shi(A)$ is the number of $\bigcirc$'s below the path.
\end{proof}

\subsubsection{$\ish=\bounce$}
The $\bounce$ statistic was discovered by Haglund in 2003 \cite{haglund03}. It provided the first combinatorial interpretation of the $q,t$-Catalan numbers of Garsia and Haiman. Haglund and Loehr \cite{haglundloehr} later extended the statistic to labeled Dyck paths $(w,\I)$ by defining $\bounce(w,\I)=\bounce(\I)$.

\begin{definition}[Haglund]
Given a Dyck path $\I$, we construct its {\sf bounce path} as follows. Begin at $(n,n)$ and travel left until we hit the path, then travel down until we hit the diagonal. Repeat these two steps until we hit $(0,0)$. Then $\bounce(\I)$ is the sum of $i$ between $1$ and $n-1$ such that the bounce path contains the diagonal point $(i,i)$.
\end{definition}

\begin{figure}
\begin{center}
\includegraphics[scale=1.3]{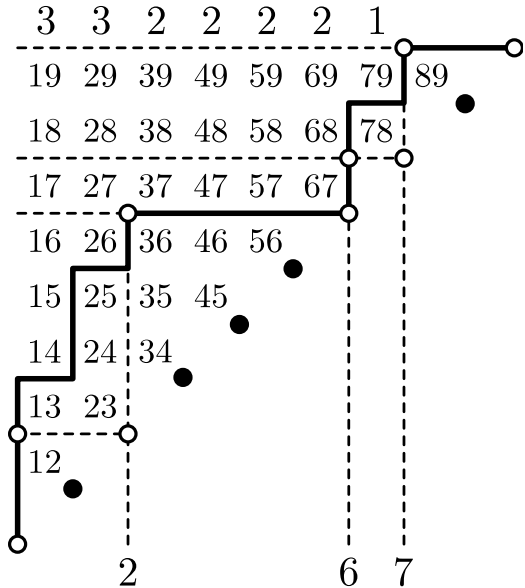}
\end{center}
\caption{The bounce path decomposition}
\label{fig:bounce_path}
\end{figure}

For example, the bounce path in Figure \ref{fig:bounce_path} is defined by the white vertices. The numbers along the bottom show that $\bounce$ for this path is $2+6+7=15$.

\begin{theorem}
Given a Shi alcove $A$ (positive or non-positive) and its corresponding labeled Dyck path $(w,\I)$ we have
\begin{equation*}
\ish(A)=\bounce(w,\I).
\end{equation*}
\end{theorem}

\begin{proof}
Suppose that $A=\w A_\circ$ where $\w$ is an affine permutation $\w\in\tilde{\S}(n)$. Suppose further that $\w=w_I \w^I$ where $w_I\in\S(n)\subseteq\tilde{\S}(n)$ is a finite permutation and $\w^I$ is a minimal coset representative. The alcove $A$ thus corresponds to a labeled Dyck path $(w_I,\I)$ and the {\em positive} alcove $\w^I A_\circ$ corresponds to the "unlabeled" Dyck path $(1,\I)$.

Recall that $\ish(A)$ is the number of hyperplanes between $\w^I A_\circ$ and $A_\circ$ of the form $e_i-e_n=a$. Given $\alpha=e_i-e_n$ the number of these hyperplanes is exactly $k_A(\alpha)$, where $k_A:\Phi^+\to\Z$ is the address of the positive Shi alcove $\w^I A_\circ$. By Theorem \ref{th:sommers}, $k_A(\alpha)=k_\I(\alpha)$ is the maximum $r$ such that $\alpha$ can be written as a sum of $r$ roots in the ideal $\I$ (i.e. above the path $\I$). For example, $\ish(A)$ is the sum of the entries in the top row of Figure \ref{fig:alcove_address}.

Now consider the bounce path of $\I$ and extend it to the left from each point that it hits $\I$. This decomposes the collection of squares above the path into ``blocks''. We have done this in Figure \ref{fig:bounce_path}; note here that there are 3 blocks. Given $\alpha\in\Phi^+$, suppose that we have $\alpha=\gamma_1+\cdots +\gamma_r$ where each $\gamma_i$ is above the path. In this case we can reorder the summands such that $\gamma_1+\cdots +\gamma_q$ is above the path for all $1\leq q\leq r$ (see \cite[Lemma 3.2]{sommers}). In other words, if $\gamma_i=e_{i_1}-e_{i_2}$, we must have $i_2=(i+1)_1$ for all $1\leq i\leq r-1$. This means there can be {\em at most one} $\gamma_i$ from each block that intersects the column containing $\alpha$.

If $\alpha=e_i-e_n$, we claim that in fact $k_\I(\alpha)$ is {\em equal to} the number of blocks that intersect the $i$th column. Indeed, set $\gamma_1=e_i-e_j$ with $j$ minimal such that $e_i-e_j\in\I$. Thus $\gamma_1$ is in the lowest block below $\alpha$. Then we travel right from $\gamma_1$, bounce off the diagonal, and travel up until we reach $\gamma_2$ such that: $\gamma_2$ is in the block above the block containing $\gamma_1$, and $\gamma_2$ is in the top row of this block. Continuing in this way, we will obtain $\alpha=\gamma_1+\cdots +\gamma_r$ such that there is one summand from each block intersecting the $i$th column. Since this was an upper bound, the claim is proved. For example, in Figure \ref{fig:bounce_path} the root $\alpha=e_1-e_9$ decomposes as $\alpha=(e_1-e_4)+(e_4-e_7)+(e_7-e_9)$, one summand from each block below $\alpha$.

Finally, $\ish(A)$ is the sum of the values $k_\I(e_i-e_n)$ for $1\leq i\leq n-1$. In other words, we sum over the number of blocks that intersect each column. This is the same as summing the number of squares in the top row of each block. Note also that there exists a block whose top row contains $j$ squares if and only if the bounce path touches the diagonal at $(j,j)$. We conclude that $\ish(A)=\bounce(\I)$.
\end{proof}

For example, the roots in the top row of Figure \ref{fig:bounce_path} have below them, respectively, 3, 3, 2, 2, 2, 2, and 1 blocks, so $\ish(A)=3+3+2+2+2+2+1=15$. On the other hand, the sum of the lengths of the top rows of the blocks is $2+6+7=15=\bounce(\I)$.

In conclusion, here is the main result of the paper.

\begin{theorem}
\label{th:main}
The bijection $A\mapsto (w,\I)$ from Shi alcoves to labeled Dyck paths sends the pair of statistics $(\binom{n}{2}-\shi,\ish)$ to the pair $(\area',\bounce)$.
\end{theorem}

\section{The Inverse Statistics}

We chose the definitions of $\shi$ and $\ish$ to emphasize their connection with the Ish hyperplane arrangement. However, we will obtain a more natural interpretation of $\ish$ when we compose it with inversion in the Weyl group. That is, let us define the following {\sf inverse statistics}.

\begin{definition}
For any affine permutation $\w\in\tilde{\S}(n)$, we define
\begin{align*}
&\shi^{-1}(\w):=\shi(\w^{-1})\quad\text{and}\\
&\ish^{-1}(\w):=\ish(\w^{-1}).
\end{align*}
\end{definition}

First let us say why we care about the inverse statistics.

\subsection{Inverse Shi Alcoves}

Let $E$ denote the set of representing alcoves for the chambers of the Shi arrangement $\Shi(n)$ (see Theorem \ref{th:sommers}). Thinking of these alcoves as elements of the affine symmetric group $\tilde{\S}(n)$ we may invert them. J.-Y. Shi showed that the set $E^{-1}$ of inverted alcoves has a remarkable shape (see \cite{shi:signtypes}).

\begin{theorem}
The inverted Shi alcoves $E^{-1}$ are precisely the alcoves inside the simplex $D^{n+1}(n)\subseteq\R_0^n$ bounded by the hyperplanes
\begin{equation*}
\{ e_i-e_{i+1}=-1 : 1\leq i\leq n-1\} \cup \{ e_1-e_n=2\},
\end{equation*}
which is congruent to the dilation $(n+1)A_\circ$ of the fundamental alcove $A_\circ$.
\end{theorem}

\begin{figure}
\begin{center}
\includegraphics[scale=1.3]{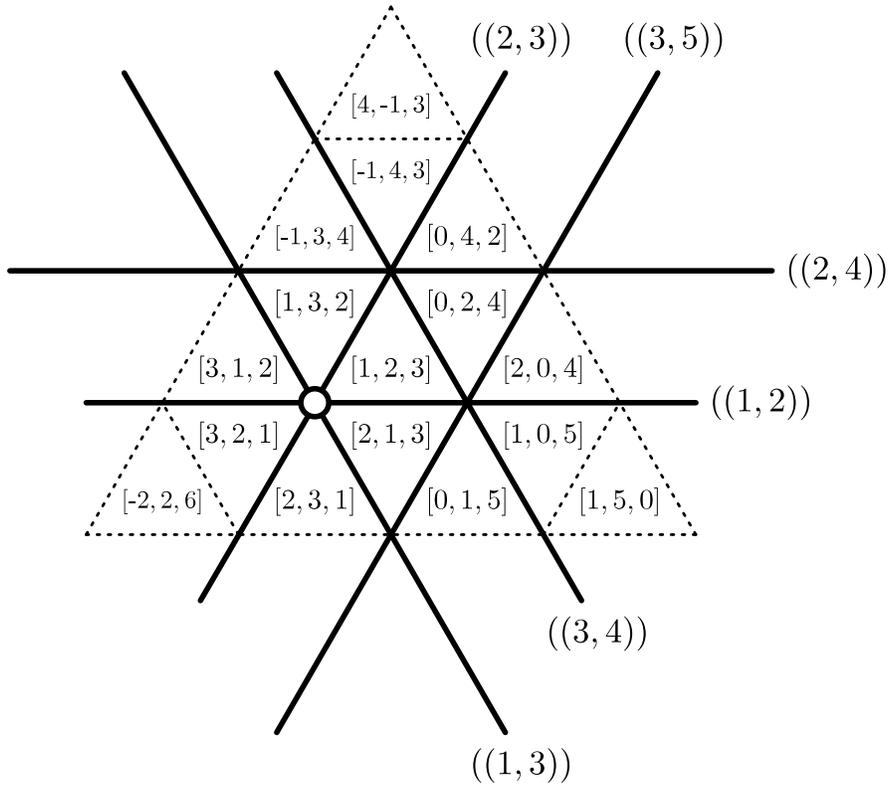}
\end{center}
\caption{The simplex $D^4(3)$ of inverted Shi alcoves}
\label{fig:shi_inverse_perms}
\end{figure}

Since the dimension of the space $\R_0^n$ is $n-1$, the simplex $D^{n+1}(n)$ contains $(n+1)^{n-1}$ alcoves. Shi concluded that his arrangement has $(n+1)^{n-1}$ chambers. Figure \ref{fig:shi_inverse_perms} displays the simplex $D^4(3)$ and the Shi arrangement in $\R_0^3$. We have labeled each alcove by the {\em inverse} of its corresponding affine permutation. Compare to Figure \ref{fig:affperm}.

Following Theorem \ref{th:main}, we assert that the joint-distribution of $\shi^{-1}$ and $\ish^{-1}$ on the simplex $D^{n+1}(n)$ is the bigraded Hilbert series of diagonal harmonic polynomials. In fact, since the shape $D^{n+1}(n)$ (Figure \ref{fig:shi_inverse_perms}) is much nicer than the distribution of Shi alcoves (Figure \ref{fig:affperm}), it seems that the inverse statistics $\shi^{-1}$ and $\ish^{-1}$ are more natural than the originals. Thus we would like to understand them directly, without reference to inversion in $\tilde{\S}(n)$.

\subsection{The Inverse $\shi$ Statistic} To do this we need to discuss the realization of the affine symmetric group $\tilde{\S}(n)$ as a semi-direct product of the finite symmetric group $\S(n)$ and the {\sf root lattice}
\begin{equation*}
Q=\{ (r(1),\ldots,r(1))\in\Z^n : r(1)+\cdots +r(n) = 0\}.
\end{equation*}
By abuse of notation, we think of $Q$ as an abelian group by associating the root $r\in Q$ with the translation $t_r:\R_0^n\to\R_0^n$ defined by $t_r(v)=v+r$. Then $\tilde{\S}(n)$ is the semi-direct product $\S(n)\ltimes Q=\{wt_r:x\in\S(n), r\in Q\}$ with multplication defined by
\begin{equation*}
(w_1t_{r_1})(w_2t_{r_2}):= (w_1w_2)(t_{w_2(r_1)}t_{r_2}) = (w_1w_2)t_{w_2(r_1)+r_2}.
\end{equation*}
Note in particular that inversion satisfies $(wt_r)^{-1}=w^{-1}t_{-w^{-1}(r)}$. 

The semi-direct product structure $\tilde{\S}(n)=\S(n)\ltimes Q$ has the following combinatorial interpretation. Recall that an affine permutation $\w:\Z\to\Z$ must satisfy $\w(k+n)=\w(k)+n$ for all $k\in\Z$, and $\w(1)+\cdots+\w(n)=\binom{n+1}{2}$. If we denote $\w\in\tilde{\S}(n)$ by the vector $\w=[\w(1),\ldots,\w(n)]$, then each affine permutation has a unique decomposition,
\begin{equation*}
[\w(1),\ldots,\w(n)] = (w(1),\ldots, w(n)) + n(r(1),\ldots,r(n))
\end{equation*}
where $w\in\S(n)$ is a finite permutation and $r=(r(1),\ldots,r(n))$ is an element of the root lattice $Q$. For example, the affine permutation $[-2,2,6]\in\tilde{\S}(3)$ decomposes as\begin{equation*}
[-2,2,6] = (1,2,3)+3(-1,0,1).
\end{equation*}
One may easily check that the map $\w=w+nr\leftrightarrow \w=wt_r$ is an isomorphism between the two structures.

We can now describe the inverse $\shi$ statistic explicitly.

\begin{theorem}For any affine permutation $\w\in\tilde{\S}(n)$ we have
\begin{equation*}
\shi^{-1}(\w)=\shi(\w).
\end{equation*}
\end{theorem}

\begin{proof}
First recall that the Shi arrangement $\Shi(n)$ consists of all the affine hyperplanes $H_{\alpha,k}$ that touch the fundamental alcove $A_\circ$. There are two of these perpendicular to each root $\alpha\in\Phi^+$; namely $H_{\alpha,0}$ and $H_{\alpha,1}$.

The inversions of $\w\in\S(n)$ are the affine hyperplanes $H$ separating the alcoves $\w A_\circ$ and $A_\circ$. These biject under the map $\w^{-1}$ to the hyperplanes $\w^{-1}H$ separating the alcoves $A_\circ$ and $\w^{-1}A_\circ$. If $\w=wt_r$, note that
\begin{equation*}
\w H_{\alpha,k} = H_{w(\alpha),k+(r,\alpha)}.
\end{equation*}
This implies that the inversions of $\w$ parallel to the root $\alpha$ biject to the inversions of $\w^{-1}$ parallel to the root $w^{-1}(\alpha)$. Finally, since every such set contains a unique Shi hyperplane (if it contains anything at all), and since the finite permutation $w^{-1}$ is a bijection on the roots $\Phi=\Phi^+\cup -\Phi^+$, we conclude that $w^{-1}$ defines a bijection from the Shi hyperplane inversions of $\w$ to the Shi hyperplane inversions of $\w^{-1}$. Hence these two sets have the same cardinality.
\end{proof}

For example, the inverse of the affine permutation $[-2,2,6]$ is $[4,2,0]$. Each of these affine permutations has $4$ inversions, among which there are $3$ Shi hyperplanes. In Section \ref{sec:nabla} we will need a more general version of this result, whose proof is the same.

\begin{theorem}
\label{th:invpart}
Given an affine permutation $\w\in\tilde{\S}(n)$, define its {\sf inversion partition} $(I_0\geq I_1\geq \ldots)$ by letting $I_k$ denote the number of affine transpositions $((i,j))$ such that $\lfloor \frac{\w(j)-\w(i)}{n}\rfloor=k$. Then $\w$ and $\w^{-1}$ have the same inversion partition.
\end{theorem}

For example, $[-2,2,6]$ and $[4,2,0]$ each have inversion partition $(3,1)$, ignoring the infinite tail of zeroes.

\subsection{The Inverse $\ish$ Statistic} Next we will compute a formula for the $\ish^{-1}$ statistic. We will find that $\ish^{-1}$ depends only on the root lattice $Q$.

To do this we need a lemma about the original $\ish$ statistic, which follows directly from Bj\"orner and Brenti \cite[Lemma 4.2]{bjornerbrenti}. The proof is instructive, so we reproduce it here.

\begin{lemma}
\label{lem:bb}
Given an affine permutation $\w\in\tilde{\S}(n)$, choose $i\in\{1,\ldots,n\}$ such that $\w(i)$ is maximum. Then
\begin{equation*}
\ish(\w)=\w(i)-n.
\end{equation*}
\end{lemma}

\begin{proof}
By definition, $\ish(\w)$ is the number of pairs $n<j$ such that $\w^I(n)>\w^I(j)$, where $\w^I$ is the affine permutation defined by taking $[\w^I(1),\ldots,\w^I(n)]$ to be the increasing rearrangement of $[\w(1),\ldots,\w(n)]$. Setting $\u=\w^I$, we have $\ish(\w)=\ish(\u)$. We wish to show that $\ish(\u)=\u(n)-n$.

So let $\u=u+nr$, where $u\in\S(n)$ is a finite permutation and $r\in Q$ is an element of the root lattice. Next fix an index $1\leq i\leq n-1$ and consider the integer $\lfloor\frac{\u(n)-\u(i)}{n}\rfloor$. Since $\u(1)<\cdots<\u(n)$, this number is always non-negative and it counts the pairs $n<i+kn$ such that $\u(n)>\u(i+kn)$. Furthermore, since $-n\leq u(n)-u(i)\leq n$, note that
\begin{equation*}
\left\lfloor\frac{\u(n)-\u(i)}{n}\right\rfloor=\left\lfloor \frac{u(n)-u(i)}{n}+r(n)-r(i)\right\rfloor
\end{equation*}
equals $r(n)-r(i)$ when $u(i)<u(n)$ and equals $r(n)-r(i)-1$ when $u(i)>u(n)$. Finally, we have
\begin{align*}
\ish(\u) &=\sum_{i=1}^{n-1} \left\lfloor\frac{\u(n)-\u(i)}{n}\right\rfloor\\
&= \left(\sum_{i=1}^{n-1} r(n)-r(i)\right)-\#\left\{1\leq i\leq n-1: u(i)>u(n)\right\}\\
&= (n-1)r(n)-\sum_{i=1}^{n-1} r(i) - (n-u(n))\\
&= (n-1)r(n)-(0-r(n))-(n-u(n))\\
&=(u(n)+nr(n))-n\\
&=\u(n)-n.
\end{align*}
\end{proof}

Somehow, everything balances to create a simple formula. From this formula we get an expression for $\ish^{-1}$.

\begin{theorem}
Given an affine permutation $\w=w+nr$, where $w\in\S(n)$ is a finite permutation and $r\in Q$ is an element of the root lattice, choose the {\em largest} index $j\in\{1,\ldots,n\}$ such that the value of $r(j)$ is a {\em minimum}. Then
\begin{equation*}
\ish^{-1}(\w)=j+n(-r(j)-1).
\end{equation*}
\end{theorem}

\begin{proof}
For $i\in\{1,\ldots,n\}$ recall that $\w^{-1}(i)=w^{-1}(i)-nr(w^{-1}(i))$. Thus the largest value of $\w^{-1}(i)$ over $i\in\{1,\ldots,n\}$ equals the largest value of $j-nr(j)$ over $j\in\{1,\ldots,n\}$. This value is achieved by the maximum $j$ such that $r(j)$ is a minimum. Lemma \ref{lem:bb} then tells us that
\begin{equation*}
\ish^{-1}(\w)=\ish(\w^{-1})=(j-nr(j))-n.
\end{equation*}
\end{proof}

\begin{figure}
\begin{center}
\includegraphics[scale=1]{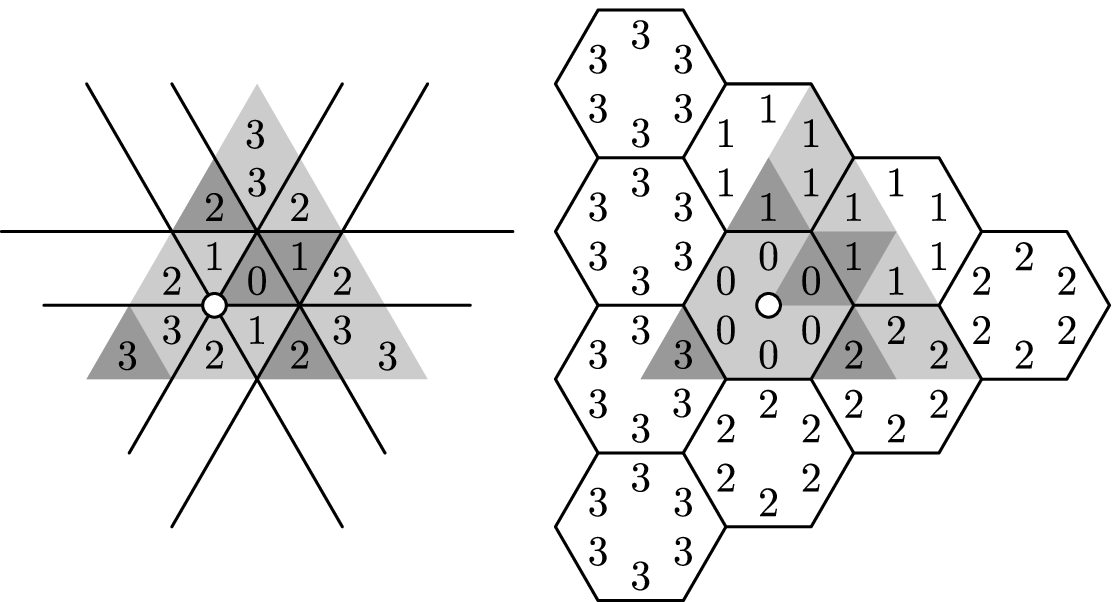}
\end{center}
\caption{The $\shi$ and $\ish^{-1}$ statistics on $D^4(3)$}
\label{fig:shi1_inverses}
\end{figure}

For example, Figure \ref{fig:shiandish} displays the $\shi$ and $\ish^{-1}$ statistics on the simplex $D^4(3)$. (The darker shaded alcoves have positive inverses.) This is the inverse of Figure \ref{fig:shi1_inverses}.

Finally, we wish to emphasize the following. The value of $\ish^{-1}(\w)$ depends only on the element $r\in Q$ of the root lattice. (This is the analogue of the fact that $\ish(\w)$ depends only on the minimal coset representative $\w^I$.) Combining this observation with Theorem \ref{th:main}, we conclude that {\em Haglund's $\bounce$ statistic is really a statistic on the root lattice of type $A$.}

\section{Powers of Nabla}
\label{sec:nabla}
In this final section we will describe several ideas for future research, roughly in order of increasing generality. Most of this depends on the {\sf nabla operator} $\nabla$ of F.~Bergeron and Garsia \cite{bergerongarsia}, which we will define first.

\subsection{The Nabla Operator} We call a formal power series in $\Q[[x_1,x_2,\ldots,]]$ a {\sf symmetric function} if it is invariant under permuting variables. Let $\Lambda=\oplus_{n\geq 0} \,\Lambda^n$ denote the ring of symmetric functions, graded by degree. Then $\Lambda^n$ is isomorphic to the space of (virtual) representations of the symmetric group $\S(n)$ over $\Q$. Under this isomorphism, the role of the irreducible representations is played by the basis of {\sf Schur functions} $s_\lambda\in\Lambda^n$, one for each partition $\lambda=(\lambda_1\geq\lambda_2\geq \cdots)$ of the integer $n=\sum_i \lambda_i$.

If we extend the field of coefficients from $\Q$ to  $\Q(q,t)$, another remarkable basis of $\Lambda^n$ is the set of {\sf modified Macdonald polynomials} $\tilde{H}_\mu$, where again $\mu=(\mu_1\geq\mu_2\geq\cdots)$ is an integer partition of $n$. Let $\nu(\mu):=\sum_{i\geq 1} (i-1)\mu_i$, and let $\mu'$ be the {\sf conjugate partition} defined by $\mu_i'=\#\{j\geq 1: \mu_j\geq i\}$. Then the Bergeron-Garsia nabla operator is the unique $\Q(q,t)$-linear map on $\Lambda^n$ defined by
\begin{equation*}
\nabla(\tilde{H}_\mu)=q^{\nu(\mu')}t^{\nu(\mu)}\tilde{H}_\mu.
\end{equation*}
That is, the modified Macdonald polynomials are a basis of eigenfunctions for $\nabla$. It turns out that many results on diagonal harmonics can be expressed elegantly in terms of $\nabla$. In particular, if $e_n=\sum_{i_1<\cdots<i_n} x_{i_1}\cdots x_{i_n}$ is the {\sf elementary symmetric function}, then $\nabla(e_n)$ is the Frobenius character of the diagonal harmonics. That is, if we replace each Schur function in $\nabla(e_n)$ by its dimension, we obtain $\DH(n;q,t)$. For details, see Haglund \cite{haglund:book}.

Now we suggest some ways to generalize our earlier results, which amount to new conjectural interpretations of the $\nabla$ operator.

\subsection{Extended Shi Arrangements} Recall that the set of reflections in the affine Weyl group $\tilde{\S}(n)$ is
\begin{equation*}
\tilde{T}=\{ ((i,j)) : \quad 1\leq i\leq n, \quad i<j \}.
\end{equation*}
The affine transposition $((i,j))$ corresponds to the hyperplane $H_{\alpha,a}$ where $\alpha=e_{\rem(j,n)}-e_{\rem(i,n)}$ and $a=\quo(j,n)$, and where {\bf remainder is taken in the set $\{1,\ldots,n\}$}. We will call $\lfloor \frac{j-i}{n}\rfloor=k$ the {\sf height} of the hyperplane; this is some measure of how far the hyperplane is from the fundamental alcove. Recall that the Shi arrangement consists of the hyperplanes of height $0$. We may now define the {\sf $m$-extended Shi arrangement}.

\begin{definition}
Let $\Shi^m(n)$ denote the set of affine transpositions $((i,j))$ with height in the set $\{0,\ldots,m-1\}$. Equivalently,
\begin{equation*}
\Shi^m(n):=\left\{ e_i-e_j=a:\,\, 1\leq i<j\leq n,\,\, a\in\{-m+1,\ldots,m\}\right\}.
\end{equation*}
\end{definition}

Athanasiadis proved \cite[Proposition 3.5]{athanasiadis} that every chamber of $\Shi^m(n)$ contains a unique alcove of minimum length. It seems true, thought we cannot find a reference, that the inverses of these representing alcoves are precisely the alcoves contained in the simplex $D^{mn+1}(n)\subseteq\R_0^n$ bounded by the hyperplanes
\begin{equation*}
\{e_i-e_{i+1}=-m-1 : 1\leq i\leq n-1\} \cup \{e_1-e_n=m+1\}.
\end{equation*}
Sommers showed in \cite[proof of Theorem 5.7]{sommers} that the simplex $D^{mn+1}(n)$ is congruent to the dilation $(mn+1)A_\circ$ of the fundamental alcove $A_\circ$. Since this occurs in the $(n-1)$-dimensional quotient space $\R_0^n$, the simplex $D^{mn+1}(n)$ consists of $(mn+1)^{n-1}$ alcoves. We would like to extend the statistics $\shi$ and $\ish^{-1}$ to these alcoves.

This turns out to be very easy to do.  The $\shi$ statistic generalizes naturally, and the $\ish$ statistic needs no generalization at all.

\begin{definition} Given an affine permutation $\w\in\tilde{\S}(n)$, let $\shi^m(\w)$ denote the number of hyperplanes of $\Shi^m(n)$ separating $\w A_\circ$ from the fundamental alcove $A_\circ$.
\end{definition}

We wish to study the joint-distribution of $\shi^m$ and $\ish$ on the minimal alcoves in the arrangement $\Shi^m(n)$. By Theorem \ref{th:invpart}, this is the same as the joint-distribution of $\shi^m$ and $\ish^{-1}$ on the alcoves of the simplex $D^{mn+1}(n)$.

\begin{conjecture}
\label{conj:pospowers}
Consider the following generating function for $\shi^m$ and $\ish^{-1}$ over alcoves in the dilated simplex $D^{mn+1}(n)$,
\begin{equation*}
\Shi^m(n;q,t):=\sum_{A\subseteq D^{mn+1}(n)} q^{\ish^{-1}(A)}t^{m\binom{n}{2}-\shi^m(A)},
\end{equation*}
and let $\Shi^m_+(n;q,t)$ denote the same sum over alcoves whose inverses are in the dominant cone $C_\circ$. We conjecture the following:
\begin{enumerate}
\item $\Shi^m(n;q,t)$ is the Hilbert series of $\nabla^m(e_n)$.
\item $q^{m\binom{n}{2}}\Shi^m(n;q,1/q)=[mn+1]_q^{n-1}$.
\item $\Shi^m_+(n;q,t)$ is the Hilbert series for the sign-isotypic component of $\nabla^m(e_n)$.
\item $q^{m\binom{n}{2}}\Shi^m_+(n;q,1/q)=\frac{1}{[n]_q}{(m+1)n \brack n-1}_q$, the $q$-Fuss-Catalan number.
\end{enumerate}
\end{conjecture}

\begin{figure}
\begin{center}
\includegraphics[scale=.65]{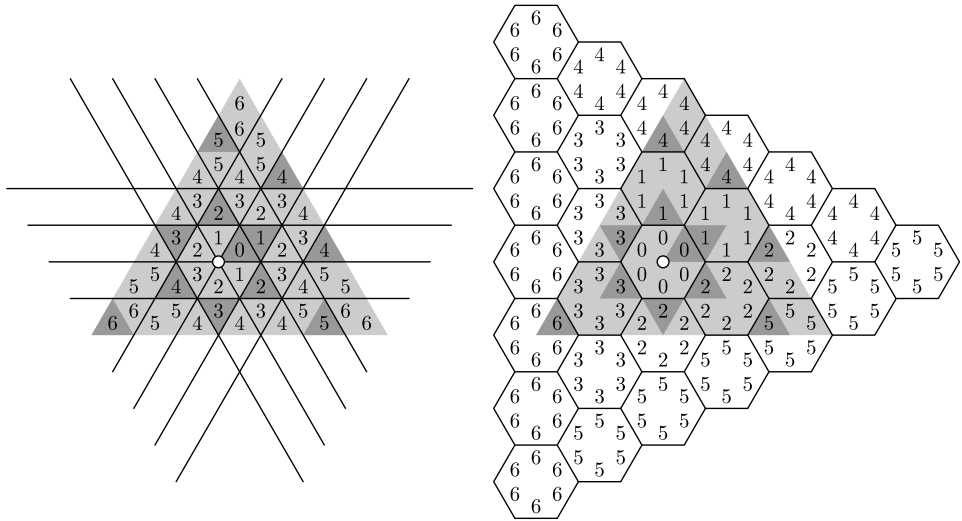}
\end{center}
\caption{The $\shi^2$ and $\ish^{-1}$ statistics on $D^7(3)$}
\label{fig:shi2_inverses}
\end{figure}

For example, let $m=2$ and $n=3$. Figure \ref{fig:shi2_inverses} displays the statistics $\shi^2$ and $\ish^{-1}$ on the alcoves of $D^7(3)$, and Table \ref{tab:shi2} displays the corresponding generating functions. One may observe that all four assertions hold in this case.

\begin{table}
\begin{center}
\begin{tabular}{cc}
\begin{tabular}{rc}
& \,\,\,\,\,\,\,\,\,\,$\ish^{-1}$ \\
\begin{tabular}{c} \vspace{.1in}\\ $\shi^2\,\,$\!\!\!\!\!\!\!\! \end{tabular} & 
\begin{tabular}{r|ccccccc}
 & 0 & 1 & 2 & 3 & 4 & 5 & 6\\
\hline
0 & 1& & & & & &\\
1 & 2& 1& & & & &\\
2 & 2& 3& 1&  & & &\\
3 & 1& 4& 3& 1 & & &\\
4 & & 3&5&3&1& &\\
5& &  1&3&4&3&1&\\
6& &&&1&2&2&1\\
\end{tabular}\\
\end{tabular}
&
\,\,\,\,\,\begin{tabular}{rc}
& \,\,\,\,\,\,\,\,\,\,$\ish^{-1}$ \\
\begin{tabular}{c} \vspace{.1in}\\ $\shi^2\,\,$\!\!\!\!\!\!\!\! \end{tabular} & 
\begin{tabular}{r|ccccccc}
 & 0 & 1 & 2 & 3 & 4 & 5 & 6\\
\hline
0 & 1& & & & & &\\
1 & & 1& & & & &\\
2 & & 1& 1&  & & &\\
3 & & & 1& 1 & & &\\
4 & & &1&1&1& &\\
5& &  &&&1&1&\\
6& &&&&&&1\\
\end{tabular}\\
\end{tabular}
\end{tabular}
\end{center}
\caption{The generating functions $\Shi^2(3:q,t)$ and $\Shi^2_+(3;q,t)$}
\label{tab:shi2}
\end{table}

Positive powers of $\nabla$ have been well-studied. We believe that a suitable extension of our main Theorem \ref{th:main} is possible, which would make our conjectures equivalent to earlier conjectures of Haiman, Loehr and Remmel (see \cite{loehrremmel}), which are based on lattice paths from $(0,0)$ to $(mn,n)$ that stay weakly above the diagonal $y=x/m$.

\subsection{Bounded Chambers} While positive powers of $\nabla$ have been investigated by several authors, to our knowledge there has been no combinatorial conjectures for {\em negative} powers of $\nabla$. In this section we will provide one.

It was shown by Edelman and Reiner \cite[Section 3]{edelmanreiner} and by Postnikov and Stanley \cite[Proposition 9.8]{postnikovstanley} that the characteristic polynomial of the $m$-extended Shi arrangement is
\begin{equation*}
\chi(\Shi^m(n),x)=(x-mn)^{n-1}.
\end{equation*}
Hence, by Zaslavsky's Theorem, $\Shi^m(n)$ has $(mn+1)^{n-1}$ chambers (as we noted above), and it has $(mn-1)^{n-1}$ bounded chambers. Athanasiadis showed that this is also an example of Ehrhart reciprocity \cite{ath:ehrhart}.

As mentioned earlier, Athanasiadis showed that each chamber of $\Shi^m(n)$ contains a unique alcove of minimum length. In the case $m=1$, Sommers showed \cite[Lemmas 5.1 and 5.2]{sommers} that, moreover, every {\em bounded} chamber of $\Shi^m(n)$ contains a unique alcove of {\em maximum} length. We believe that this is true for general $m\geq 1$, and furthermore we believe that the inverses of these alcoves are precisely the alcoves contained in the simplex $D^{mn-1}(n)\subseteq\R_0^n$ bounded by the hyperplanes
\begin{equation*}
\{e_i-e_{i+1}=m : 1\leq i\leq n-1\} \cup \{e_1-e_n=-m+1\}.
\end{equation*}
As with $D^{mn+1}(n)$ above, Sommers has shown that $D^{mn-1}(n)$ is congruent to the dilation $(mn-1)A_\circ$ of the fundamental alcove, which implies that $D^{mn-1}(n)$ contains $(mn-1)^{n-1}$ alcoves. We wish to study the statistics $\shi^m$ and $\ish^{-1}$ on these alcoves.

\begin{conjecture}
\label{conj:negpowers}
Consider the following generating function for $\shi^m$ and $\ish^{-1}$ over alcoves in the dilated simplex $D^{mn-1}(n)$,
\begin{equation*}
\Shi^{-m}(n;q,t):=\sum_{A\subseteq D^{mn-1}(n)} q^{\ish^{-1}(A)}t^{(mn-2)(n-1)/2-\shi^m(A)},
\end{equation*}
and let $\Shi_+^{-m}(n;q,t)$ denote the same sum over alcoves whose inverses are in the dominant cone $C_\circ$. We conjecture the following.
\begin{enumerate}
\item $(-1)^{n-1}\Shi^{-m}(n;1/q,1/t)/q^{n-1}t^{n-1}$ is the Hilbert series of $\nabla^{-m}(e_n)$.
\item $q^{(mn-2)(n-1)/2}\Shi^{-m}(n;q,1/q)=[mn-1]_q^{n-1}$.
\item $(-1)^{n-1}\Shi_+^{-m}(n;1/q,1/t)/q^{n-1}t^{n-1}$ is the Hilbert series for the sign-isotypic component of $\nabla^{-m}(e_n)$.
\item $q^{(mn-2)(n-1)/2}\Shi_+^{-m}(n;q,1/q)=\frac{1}{[n]_q}{(m+1)n-2 \brack n-1}_q$.
\end{enumerate}
\end{conjecture}

\begin{figure}
\begin{center}
\includegraphics[scale=.65]{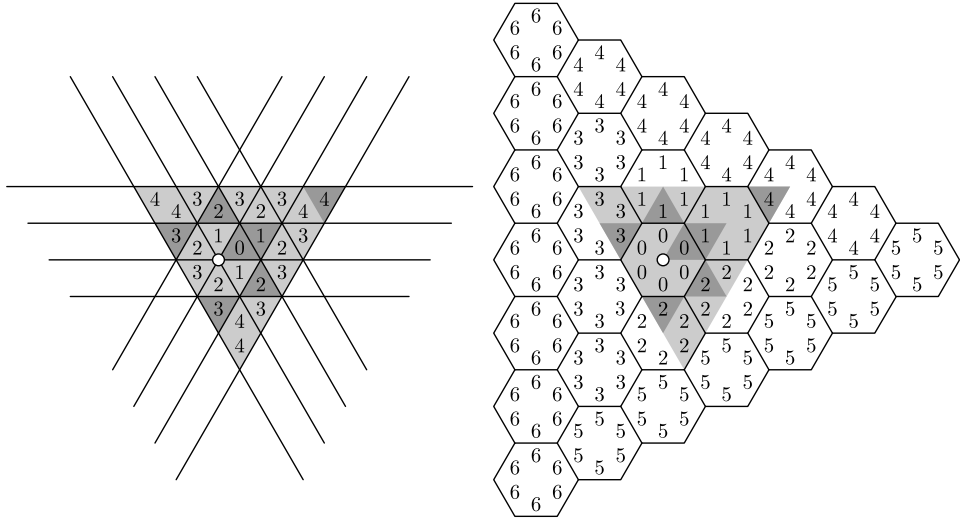}
\end{center}
\caption{The $\shi^2$ and $\ish^{-1}$ statistics on $D^5(3)$}
\label{fig:shi2_bounded_inverses}
\end{figure}

For example, Figure \ref{fig:shi2_bounded_inverses} displays the statistics $\shi^2$ and $\ish^{-1}$ on the simplex $D^{5}(3)$, and Table \ref{tab:shi2_bounded} displays the corresponding generating functions.  One may observe that all four assertions hold for this data.

\begin{table}
\begin{center}
\begin{tabular}{cc}
\begin{tabular}{rc}
& \,\,\,\,\,\,\,\,\,\,$\ish^{-1}$ \\
\begin{tabular}{c} \vspace{.1in}\\ $\shi^2\,\,$\!\!\!\!\!\!\!\! \end{tabular} & 
\begin{tabular}{r|ccccc}
 & 0 & 1 & 2 & 3 & 4\\
\hline
0 & 1& & & &\\
1 & 2& 1& & & \\
2 & 2& 3& 1& &\\
3 & 1& 4& 3& 1 &\\
4 &  & 1&2&2&1
\end{tabular}\\
\end{tabular}
&
\,\,\,\,
\begin{tabular}{rc}
& \,\,\,\,\,\,\,\,\,\,$\ish^{-1}$ \\
\begin{tabular}{c} \vspace{.1in}\\ $\shi^2\,\,$\!\!\!\!\!\!\!\! \end{tabular} & 
\begin{tabular}{r|ccccc}
 & 0 & 1 & 2 & 3 & 4\\
\hline
0 & 1& & & &\\
1 & & 1& & & \\
2 & & 1& 1& &\\
3 & & & 1& 1 &\\
4 &  & &&&1
\end{tabular}\\
\end{tabular}
\end{tabular}
\end{center}
\caption{The generating functions $\Shi^{-2}(3:q,t)$ and $\Shi^{-2}_+(3;q,t)$}
\label{tab:shi2_bounded}
\end{table}

Combining Conjectures \ref{conj:pospowers} and \ref{conj:negpowers}, we obtain a conjectural combinatorial interpretation {\em for all integral powers of the nabla operator acting on $e_n$}. We wonder whether Athanasiadis' result \cite{ath:ehrhart} on Ehrhart reciprocity for Shi arrangements may reflect some sort of reciprocity theorem for the nabla operator.

\subsection{Interpolation} Since the forms of Conjectures \ref{conj:pospowers} and \ref{conj:negpowers} are so similar, one may ask whether there is a more general form encompassing them both. In this case we do not have a concrete conjecture, but we we will suggest some ideas.

The simplices $D^{mn+1}(n)$ and $D^{mn-1}(n)$ are both special cases of the following construction of Sommers.

Recall that the root system of type $A_{n-1}$ is defined by $\Phi=\{e_i-e_j : 1\leq i,j\leq n \}$, and the basis of simple roots is $\Delta=\{e_i-e_{i+1}: 1\leq i\leq n-1\}$.  Given a root $\alpha\in\Phi$, let $b$ denote the sum of its coefficients in the simple root basis; we say that $b$ is the {\sf height} of the root $\alpha$. Let $\Phi_b\subseteq\Phi$ denote the set of roots of height $b$. Finally, let $p=an+b$ be any integer coprime to $n$, with $1\leq b\leq n-1$, and let $D^p(n)$ be the region containing the origin and bounded by the hyperplanes
\begin{equation*}
\{ H_{\alpha,a} : \alpha\in\Phi_{-b}\} \cup \{ H_{\alpha,a+1} : \alpha\in\Phi_{n-b}\}.
\end{equation*}
As with $D^{mn+1}(n)$ and $D^{mn-1}(n)$, Sommers showed that for any $p$ coprime to $n$, $D^p(n)$ is congruent to the dilation $pA_\circ$ of the fundamental alcove; that it contains $p^{n-1}$ alcoves; and that it contains $\frac{1}{p+n}\binom{p+n}{n}$ alcoves whose inverses are in the dominant cone. We suggest the following:

\begin{problem}
Define a statistic $\stat$ on the alcoves of $D^p(n)$. Consider the generating function
\begin{equation*}
F(p,n;q,t):=\sum_{A\subseteq D^p(n)} q^{\ish^{-1}(A)} t^{(p-1)(n-1)/2 -\stat(A)},
\end{equation*}
and let $F_+(p,n;q,t)$ denote the same sum over alcoves whose inverses lie in the dominant cone $C_\circ$. These generating functions should satisfy
\begin{enumerate}
\item $F(p,n;q,t)=F(p,n;t,q)$.
\item $q^{(p-1)(n-1)/2}F(p,n;q,1/q)=[p]_q^{n-1}$.
\item $q^{(p-1)(n-1)/2}F(p,n;q,1/q)=\frac{1}{[p+n]_q}{p+n \brack n}_q$.
\end{enumerate}
\end{problem}

Note that $\ish^{-1}$ does not need to be modified. It is the $\shi$ statistic that is difficult to define in general. We note that the smallest mystery case is $p=2$ and $n=5$, which corresponds to the $4$-dimensional simplex $D^2(5)$ in $\R_0^5$ bounded by the hyperplanes
\begin{equation*}
e_1-e_3=e_2-e_4=e_3-e_5=0 \quad\text{ and }\quad e_1-e_4=e_2-e_5=1.
\end{equation*}
This simplex contains $2^4=16$ alcoves, corresponding to the affine permutations
\begin{equation*}
\begin{tabular}{rrrr}
$[-1,2,5,3,6],$ & $[0,3,2,4,6],$ & $[1,2,4,3,5],$ & $[2,1,3,4,5],$ \\
$[0,2,3,4,6],$ & $[2,0,3,6,4],$ & $[1,3,2,4,5],$ & $[2,1,3,5,4],$ \\
$[0,2,4,3,6],$ & $[1,2,3,4,5],$ & $[1,3,2,5,4],$ & $[2,1,4,3,5],$ \\
$[0,3,1,4,6],$ & $[1,2,3,5,4],$ & $[1,4,2,5,3],$ & $[3,1,4,2,5]$.
\end{tabular}
\end{equation*}
Of these, only $\frac{1}{2+5}\binom{2+5}{5}=3$ have inverses in the dominant cone --- namely, $[1,2,3,4,5]$, $[0,2,3,4,6]$ and $[2,0,3,6,4]$. The distribution of $\ish^{-1}$ over the former $16$ is $\sum q^{\ish^{-1}(A)}=10+5q+q^2$ and the distribution of $\ish^{-1}$ over the latter $3$ is $\sum q^{\ish^{-1}(A)}=1+q+q^2$. We do not know what the analogue of $\shi$ is in this case, but we have checked that it cannot simply be the number of inversions coming from some special set of hyperplanes.

\subsection{Other Types} In this paper we have focused on the affine Weyl group of type $\tilde{A}_n$, which is the group $\tilde{\S}(n)$ of affine permutations. However, we have tried to use language throughout that is general to all affine Weyl groups. Certainly, the combinatorics of Shi arrangements is completely general. Also, Sommers' simplex$D^p(h)$ is defined in general for any integer $p$ coprime to the Coxeter number $h$.

Haiman observed that the most obvious generalization of the ring of harmonic polynomials to other types is ``too large'' (see \cite[Section 7]{Hai94}), and he conjectured that some suitable quotient should be considered instead. Using rational Cherednik algebras, Gordon \cite{gordon} was able to construct such a quotient. Gordon and Griffeth \cite{gordongriffeth} have now observed that this module does have a suitable bigrading and it satisfies many of the desired combinatorial properties. Gordon-Griffeth \cite{gordongriffeth}  and Stump \cite{stump} have both defined $q,t$-Catalan numbers in general type (even in complex types), however their numbers disagree in the non-well-generated complex types. This is an active area, and we wish to emphasize: {\bf as of this writing, there is no known combinatorial interpretation for these objects beyond type $A$}.

We suggest that the $\shi$ statistic on the simplex $D^{h+1}(h)$ is a good place to start. The next step is to define an analogue of the $\ish^{-1}$ statistic. Unfortunately, we have checked that in type $B_2$ it cannot simply be a statistic on the root lattice.

\section{Acknowledgements}

The author thanks Christos Athanasiadis, Steve Griffeth, Nick Loehr, Eric Sommers, and Greg Warrington for helpful discussions. Extra thanks are due to Eric Sommers for suggesting the proof of Theorem \ref{th:invpart} and to Greg Warrington for providing \texttt{Maple} code for working with $\nabla$. The original idea for this work was inspired by discussions with Steve Griffeth and by the paper \cite{fishelvazirani} of Fishel and Vazirani.

\end{document}